# A "nearly parametric" solution to Selective Harmonic Elimination PWM

Bao-Xin Shang[1], Shu-Gong Zhang[2], Na Lei[3]*, Jing-Yi Chen[4]

**Abstract**

Selective Harmonic Elimination Pulse Width Modulation (SHEPWM) is an important technique to solve PWM problems, which control the output voltage of an inverter via selecting appropriate switching angles. Based on the Rational Univariate Representation (RUR) theory for solving polynomial systems, the paper presents an algorithm to compute a "nearly parametric" solution to a SHEPWM problem. When the number of switching angles *N* is fixed, a "nearly parametric" solution can be considered as functions of the modulation index *m*. So we can adapt the amplitude of the output voltage with the same source voltage by changing the modulation index. When *m* is given as a specific value, complete solutions to the SHEPWM problem can be obtained easily using univariate polynomial solving. Compared with other methods, *m* is considered as a symbolic parameter for the first time, and this can help avoid totally restarting when *m* changes. The average time for computing complete solutions associated to 460 modulation indexes based on a "nearly parametric" solution when *N*=5 is 0.0284s, so the algorithm is practical. Three groups of switching angles associated to *N*=5, *m*=0.75 is simulated in MATLAB, and it verifies the algorithm's correctness.

## I. INTRODUCTION

Multilevel inverters have lots of applications in the field of industry [1]. Pulse Width Modulation (PWM) is of great importance in medium-voltage high-power inverters, which controls the output voltage of an inverter through selecting appropriate switching angles. [1]-[2] review topologies and modulation techniques of multilevel inverters. SHEPWM is a widely-used method to solve PWM problems, and it has advantages such as low switching loss, linearly controlling fundamental voltage component, and specifying the orders of eliminated harmonics. The method is firstly introduced in [3], which use the Newton-Raphson iteration method to solve associated nonlinear equations. Abdul Moeed Amjad and Zainal Salam summarize a variety of methods to solve SHEPWM problems in [4]. With the number of switching angles *N* fixed, the difficulty of a SHEPWM problem lies in how to solve the associated transcendental system. There mainly exist three difficulties. The first one is how to get a solution given a specific modulation index *m*. In this case, a group of good initial values is important. Initial values selection techniques are considered in [5]-[7]. In addition, a great number of optimization algorithms such as genetic method, bee method, and warm particle method etc are used in [8]-[11], and they can find more than one group of switching angles and allow rough initial values. The second one is how to get complete solutions given a specific modulation index *m*. In 2004, John Chiasson etc present a method to get complete solutions based on resultant theory in symbolic computation[12]. In 2007, Qunjing Wang etc improve the result using symmetric polynomial theory[13]. The last difficulty is when the modulation index *m* changes how to recompute related complete solutions quickly. If this comes true, output voltage can be adapted more easily through modifying *m* with source voltage unchanged. To this end, all above methods have to restart totally. To fix this problem, we use the RUR method to solve parametric zero-dimensional polynomial systems, and firstly consider utilizing algebraic interpolation techniques to accelerate the computation. We implements the algorithm in Maple 18. When *m* changes, recomputation starts from a "nearly parametric" solution, and it only involves little times of univariate polynomial solving. More information on theories and implementation of RUR is referred to [14]-[15].

The structure of the paper is as follows. In section II, we describe how to transform the transcendental equations associated to SHEPWM problems into polynomial equations in elementary symmetric polynomials. We also implements a fast conversion algorithm in Maple 18 to accelerate the transforming process, which transforms power sum symmetric polynomials in $x_i$ to polynomials in $s_i$, $i = 1, \ldots, N$, where $s_i$'s are elementary symmetric polynomials in $x_i$'s. Section III contains theories on the RUR method, and we use algebraic interpolation to accelerate the computation of a RUR of a parametric zero-dimensional polynomial system. This help us get a "nearly parametric" solution to a SHEPWM problem. Section IV depicts the steps to solve a SHEPWM problem with the number of switching angles *N* fixed. It mainly contains two parts. The first part focuses on preprocessing, and gets a "nearly parametric" solution to a SHEPWM problem, in which the modulation index *m* is considered as a parameter. The other part shows how to get complete solutions when *m* is given as a specific value. Section V lists numeric experiments and simulation results. At last, the algorithm described in section IV is illuminated by a detailed

[1] Email: shbxin@163.com. School of Mathematics, Jilin University, Changchun 130012, PR China; College of Science, Northeast Dianli University, Jilin 132012, PR China
[2] School of Mathematics, Jilin University, Changchun 130012, PR China
[3] Email: leina@jlu.edu.cn. School of Mathematics, Jilin University, Changchun 130012, PR China
[4] Electrical Engineering College, Northeast Dianli University, Jilin 132012, PR China

example when N=3 in appendix A. Appendix B lists some facts about symmetric polynomials. In appendix C, the proof of Theorem 1 is given.

## II. MODEL DESCRIPTION

We mainly consider half-wave symmetry SHEPWM problems of three-level voltage inverters. Consider the unipolar switching scheme, and the bipolar case is similar. Fig. 1 shows one phase voltage waveform when the number of switching angles $N$ is even; Fig. 2 shows that when $N$ is odd. The aim of SHEPWM problems is to determine $N$ switching angles $0 \leq \alpha_1 < \cdots < \alpha_N < \frac{\pi}{2}$ to eliminate $N-1$ harmonics in the output voltage.

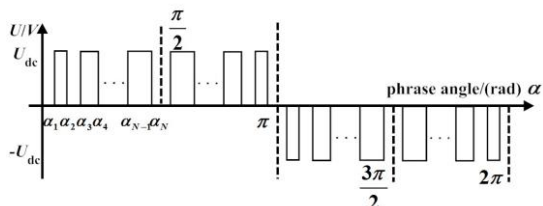

Fig. 1. One phase voltage waveform ($N$ is even).

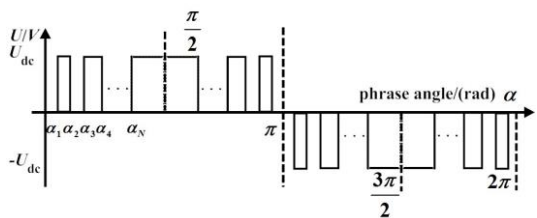

Fig. 2. One phase voltage waveform ($N$ is odd).

Without loss of generality, we can assume the amplitude of $U(\alpha)$ is 1. Apply Fourier analysis to $U(\alpha)$ and get

$$U(\alpha) = \sum_{k=1}^{\infty} b_k \sin k\alpha.$$

Due to the symmetry property, the Fourier coefficients

$$b_k = \frac{4}{\pi} \int_0^{\pi/2} U(\alpha) \sin k\alpha \, d\alpha.$$

Whenever $N$ is even or odd, $b_k = \frac{4}{k\pi} \sum_{j=1}^{N}(-1)^{j+1} \cos k\alpha_j$, $k = 1, 3, 5, \ldots$. We consider how to eliminate 5th, 7th, 11th, 13th, 17th, 19th,…harmonics. Since there are $N$ switching angles, we can eliminate $N-1$ selective harmonics. Set $b_1 = \frac{4}{\pi} m$, where $0 < m < 1$ is called the *modulation index*. $m$ is a measure of the magnitude of the fundamental output voltage. To eliminate $i^{th}$ harmonic, set $b_i = 0, i = 5, 7, 11, 13, \ldots$. Then we get equations

$$\begin{cases} \frac{4}{\pi}\left(\sum_{i=1}^{N}(-1)^{i+1}\cos\alpha_i - m\right) = 0, \\ b_{6M\pm 1} = \frac{4}{(6M\pm 1)\pi}\sum_{i=1}^{N}(-1)^{i+1}\cos(6M\pm 1)\alpha_i = 0,\ M=1,\ldots,\left\lfloor\frac{N}{2}\right\rfloor - 1, \\ b_{6\left\lfloor\frac{N}{2}\right\rfloor\pm 1} = \frac{4}{(6\left\lfloor\frac{N}{2}\right\rfloor\pm 1)\pi}\sum_{i=1}^{N}(-1)^{i+1}\cos(6\left\lfloor\frac{N}{2}\right\rfloor\pm 1)\alpha_i = 0. \end{cases}$$

(1)

If $N$ is odd, "$\pm$" is used in the last expression; if $N$ is even, "$-$" is used.

Note that $\cos k\alpha, k = 2, 3, \ldots$ can be represented as a polynomial in $\cos\alpha$. Set $x_i = (-1)^{i+1}\cos\alpha_i$, then equations (1) can be converted into symmetric polynomial equations in $x_1, \ldots, x_N$ [12]-[13]. And we denote them orderly by

$$g_1(x_1, \ldots, x_N) = \sum_{i=1}^{N} x_i - m = 0,$$
$$g_2(x_1, \ldots, x_N) = 0, \qquad (2)$$
$$\ldots$$
$$g_N(x_1, \ldots, x_N) = 0.$$

The reason why $x_i = (-1)^{i+1}\cos\alpha_i$ rather than $x_i = \cos\alpha_i$ are used is to get symmetric polynomials in $x_1, \ldots, x_N$. Now, the problem is how to solve polynomial equations (2). Considering the symmetric property of the equations (2), it can be converted into equations in $s_1, \ldots, s_N$, which denote the elementary symmetric polynomials in $x_1, \ldots, x_N$ [13]. Denote them orderly by

$$f_1(s_1, \ldots, s_N) = s_1 - m = 0, \ldots, f_N(s_1, \ldots, s_N) = 0. \quad (3)$$

Both Mathematica 9 and Maple 18 have build-in procedures to implement above conversion process. But for this problem, both of them have efficiency problems. Note that polynomials in equations (2) are sums of power sum symmetric polynomials, so the explicit conversion formula from power sum symmetric polynomials to elementary symmetric polynomials is much more efficient. We implement the conversion algorithm in Maple 18, and it consumes much less CPU time and memory than build-in procedures in Mathematica 9 and Maple 18.

For our problem, there are two main advantages to do the conversion. The first is to reduce the number of variables by 1. The second is to reduce the number of zeros of polynomial equations to be solved significantly. They will be illuminated by the detailed example in Appendix A. Next, we will give a theorem based on which useless computation can be reduced.

If a monic univariate polynomial of degree $n$ has $n$ distinct real zeros satisfying $x_1 > -x_2 > \cdots > (-1)^{n+1}x_n > 0$, then the signs of its coefficients are determined. It can be stated as follows.

**Theorem 1** Suppose that $x_1 > -x_2 > \cdots > (-1)^{n+1}x_n > 0$. Let $s_1 = \sum_{i=1}^{n} x_i, \ldots, s_r = \sum_{1 \leq i_1 < i_2 < \cdots < i_r \leq n} x_{i_1} x_{i_2} \cdots x_{i_r}, \ldots, s_n = x_1 \cdots x_n$, $1 \leq r \leq n$. Then $s_{4d+1} > 0, s_{4d+2} < 0, s_{4d+3} < 0, s_{4d+4} > 0$, $d = 0, 1, \ldots$.

Theorem 1 can help us avoid useless computation, and it is explained in the second part of the algorithm in section IV. A proof of Theorem 1 is listed in Appendix C.

The following section will depict the tool used to solve polynomial equations (3) in $s_1, \ldots, s_N$, which contains *m* as a parameter.

## III. RATIONAL UNIVARIATE REPRESENTATION FOR SOLVING POLYNOMIAL SYSTEMS

### A. RURs for zero-dimensional polynomial systems

Rational Univariate Representation (RUR) is an efficient method for solving zero-dimensional polynomial systems which only have finite zeros. In 1999, Rouillier details the method in [14], and gives an efficient way to compute such a representation. In essence, it constructs a one-to-one map between zeros of the original system and those of a univariate polynomial $\chi_t(T)$. Then the zeros of the original system can be represented by those of $\chi_t(T)$. By now, this method provides closer solutions to "parametric" ones for a general zero-dimensional polynomial than any other symbolic or numeric methods. More information about RUR is referred to [14, 15, 16, 17].

Let $\mathbb{C}$ be the field of complex numbers and $\mathbb{Q}$ the field of rational numbers. Let $f_1, \ldots, f_s \in \mathbb{Q}[x_1, \ldots, x_n]$, which have finite zeros in $\mathbb{C}^n$. Let *t* be a $\mathbb{Q}$–linear combination of $x_1, \ldots, x_n$. *t* is called a separating element if *t* has distinct values at zeros of $f_1, \ldots, f_s$. Then, the RUR of $f_1, \ldots, f_s$ with respect to *t* is

$$\{\chi_t(T), g_t(1,T), g_t(x_1,T), \ldots, g_t(x_n,T)\},$$

where $\chi_t(T), g_t(1,T), g_t(x_1,T), \ldots, g_t(x_n,T) \in \mathbb{Q}[T]$. The RUR is unique up to *t*. For convenience, we also call it a RUR of $f_1 = 0, \ldots, f_s = 0$. Then the zeros of $f_1, \ldots, f_s$ can be expressed as

$$\left\{ \left( \frac{g_t(x_1,T_0)}{g_t(1,T_0)}, \ldots, \frac{g_t(x_n,T_0)}{g_t(1,T_0)} \right) \middle| \chi_t(T_0) = 0, T_0 \in \mathbb{C} \right\}.$$

Moreover, there exists a bijection between the rational zeros (real zeros, complex zeros) of $\chi_t(T)$ and those of $f_1, \ldots, f_s$.

In Maple 18, the built-in procedure RationalUnivariateRepresentation is used to compute a RUR of a zero-dimensional polynomial system. For our problem, due to the changeability of *m*, directly using RationalUnivariateRepresentation is impossible. So we consider RURs of parametric zero-dimensional polynomial systems.

### B. RURs for parametric zero-dimensional polynomial systems

Since polynomial systems in the paper involve only one parameter, we limit the number of parameters to 1, and denote it by *u*. Let $\mathbb{Q}[u]$ be the set of all polynomials in *u* with coefficients in $\mathbb{Q}$, $\mathbb{Q}(u)$ the set of all rational functions with their numerators and denominators in $\mathbb{Q}[u]$. Let $f_1, \ldots, f_s \in \mathbb{Q}[u][x_1, \ldots, x_n]$. If $f_1 = 0, \ldots, f_s = 0$ has only finite zeros in the algebraic closure of $\mathbb{Q}(u)$, then $f_1, \ldots, f_s$ is called a parametric zero-dimensional polynomial system with *u* a parameter. Then the RUR of $f_1, \ldots, f_s$ with respect to *u*, *t* is represented as

$$\{\chi_{t,u}(T), g_{t,u}(1,T), g_{t,u}(x_1,T), \ldots, g_{t,u}(x_n,T)\},$$

where $\chi_{t,u}(T), g_{t,u}(1,T), g_{t,u}(x_1,T), \ldots, g_{t,u}(x_n,T) \in \mathbb{Q}(u)[T]$. $\chi_{t,u}(T)$ is monic, and *t* is a linear combination of $x_1, \ldots, x_n$, which is a separating element of $f_1, \ldots, f_s$. The computing of separating elements is referred to [15, 16]. For $u_0 \in \mathbb{C}$, if $u_0$ satisfies:

(1) $\chi_{t,u}(T)|_{u=u_0}, g_{t,u}(1,T)|_{u=u_0}, g_{t,u}(x_1,T)|_{u=u_0}, \ldots, g_{t,u}(x_n,T)|_{u=u_0}$ are well-defined.

(2) $\overline{\chi_{t,u}(T)}\big|_{u=u_0}$ is square-free, where $\overline{\chi_{t,u}(T)}$ denotes the square-free part of $\chi_{t,u}(T)$.

Then

$$\left\{ \left( \frac{g_t(x_1,T_0)}{g_t(1,T_0)}\bigg|_{u=u_0}, \ldots, \frac{g_t(x_n,T_0)}{g_t(1,T_0)}\bigg|_{u=u_0} \right) \middle| \chi_t(T_0)|_{u=u_0} = 0, T_0 \in \mathbb{C} \right\}$$
(4)

is the zeros set of $f_1|_{u=u_0}, \ldots, f_s|_{u=u_0}$. More discussions about RURs of parametric zero-dimensional polynomial systems are referred to [15].

As far as I know, there is no built-in procedures in Maple 18 or Mathematica 9 to compute a RUR of a parametric zero-dimensional polynomial system. The algorithm in [14, 15] can be used to solve this problem, but its computation involves rational functions, so it consumes too much CPU time and RAM. To fix this problem, note that the coefficients with respect to *T* of polynomials in RURs of parametric zero-dimensional polynomial systems are rational functions in *u*, so we can use algebraic interpolation theory to improve the algorithm in [15], and the improved algorithm computes RURs with operations all over the field of rational numbers. We implements the improved algorithm in Maple 18.

A RUR of equations (3) with respect to *m, t* is called a *"nearly parametric"* solution to the SHEPWM problem, which is only related to the number of switching angles. When *m* changes, computation only need to restart from the "nearly parametric" solution, and this avoids much useless computation. In the next section, we show in detail the steps to compute "nearly parametric" solutions to SHEPWM problems and how to get complete solutions based on a "nearly parametric" solution.

## IV. SOLVING SHEPWM PROBLEMS

Fix the number of switching angles *N* and consider the modulation index *m* as a parameter. The algorithm to solve SHEPWM problems is split into two parts. The first part mainly includes preprocessing, and its target is to get a "nearly parametric" solution to a SHEPWM problem, which is a RUR of equations (3). The second part is to get complete solutions to equations (1) based on the "nearly parametric" solution when *m* is given as a specific value. In this process, only univariate polynomial solving is involved. So if *m* changes, we only need to restart the computation from the "nearly parametric" solution, and computation cost is cut down.

The first part of the algorithm is described as follows.

(1) Given the number of switching angles *N*, get

transcendental equations (1).

(2) Let $x_i = (-1)^{i+1}\cos\alpha_i, i=1,\ldots,N$ and convert equations (1) into polynomial equations in $x_1,\ldots,x_N$. Denote them orderly by $g_1 = \sum_{i=1}^N x_i - m = 0, g_2 = 0, \ldots, g_N = 0$. Then using elementary symmetric polynomials $s_1,\ldots,s_N \in \mathbb{Q}[x_1,\ldots,x_n]$, translate $g_1=0, g_2=0,\ldots,g_N=0$ into polynomial equations in $s_1,\ldots,s_N$, and we get

$$\begin{cases} f_1(s_1,\ldots,s_N) = s_1 - m = 0, \\ f_2(s_1,\ldots,s_N) = 0, \\ \cdots \\ f_N(s_1,\ldots,s_N) = 0. \end{cases} \quad (5)$$

(3) According to the physical significance of $m$, $s_1$ or $m$ can be considered as a parameter of the last $N-1$ equations of (5). Using the theory of RUR, compute a RUR of equations (5) with respect to $m$, and we get a "nearly parametric" solution to the SHEPWM problem

$$\chi_{t,m}(T), g_{t,m}(1,T), g_{t,m}(s_2,T),\ldots, g_{t,m}(s_N,T) \quad (6)$$

When the modulation index $m$ is fixed, we can use formula (6) to get $s_2,\ldots,s_N$. As a result of Theorem 1, only $s_{4d+1} > 0, s_{4d+2} < 0, s_{4d+3} < 0, s_{4d+4} > 0, d \ge 0$ has a possibility to derive an effective group of switching angles.

For a fixed modulation index $m_0$, the second part of the algorithm is depicted as follows.

(1) Substitute $m = m_0$ into formula (6), and get

$$\chi_{t,m}^{(m_0)}(T), g_{t,m}^{(m_0)}(1,T), g_{t,m}^{(m_0)}(s_2,T),\ldots, g_{t,m}^{(m_0)}(s_N,T). \quad (7)$$

(2) Solve $\chi_{t,m}^{(m_0)}(T) = 0$ and denote its roots by $T_1,\ldots,T_l$.

(3) Substitute $T = T_i, i=1,\ldots,l$ into $s_2 = \dfrac{g_{t,m}^{(m_0)}(s_2,T)}{g_{t,m}^{(m_0)}(1,T)}$, $\ldots, s_N = \dfrac{g_{t,m}^{(m_0)}(s_N,T)}{g_{t,m}^{(m_0)}(1,T)}$, and denote the results by $S_{i'} = (m_0, s_{i',2},\ldots,s_{i',N}), i'=1,\ldots,l', l' \le l;$ $s_{i',2},\ldots,s_{i',N}$ satisfy $s_{i',4d+1} > 0, s_{i',4d+2} < 0, s_{i',4d+3} < 0, s_{i',4d+4} > 0, d \ge 0.$ Note that we abandon the $s_i's$ whose signs don't conform to the conclusion of Theorem 1.

(4) According to Vieta's formulas, for $S_{i'}$, solve $x^N - m_0 x^{N-1} + s_{i',2} x^{N-2} + \cdots + (-1)^N s_{i',N} = 0$, then we can get $x_1^{(i')},\ldots,x_N^{(i')}$.

(5) If there exists a permutation $x_{j_1}^{(i')},\ldots,x_{j_N}^{(i')}$ of $x_1^{(i')},\ldots,x_N^{(i')}$ satisfying $1 \ge x_{j_1}^{(i')} > \cdots > (-1)^{N+1} x_{j_N}^{(i')} > 0$, we get a group of switching angles $\alpha_k = \arccos(-1)^{1+k} x_{j_k}^{(i')}$.

Only steps (2) and (4) involves univariate polynomial solving, so the second part of the algorithm performs at most $l+1$ times univariate polynomial solving.

We give a detailed example in Appendix A to illuminate how the algorithm works.

## V. NUMERIC EXPERIMENTS AND SIMULATION

Based on a "nearly parametric" solution with $N=5$, we compute complete solutions to SHEPWM problems associated to $m = i/500, i=1,\ldots,460$ in Maple 18. And it averagely consumes 0.0284s on a laptop with an Intel Dual Core P8400 (2 × 2.26GHz) and 4G RAM.

We draw the switching angle trajectories in Fig. 3 when $N=5$. It consists of 1035 groups of switching angles associated to $m=i/500, i=1,\ldots,460$. The points with the same $m$ and color form a group of switching angles. According to the experiment, there exist two groups of switching angles when $0 < m \le 0.478$, $0.516 \le m \le 0.528$, and $0.786 \le m \le 0.918$; three groups of switching angles when $0.479 \le m \le 0.487$ and $0.529 \le m \le 0.785$; only one group of switching angles when $0.488 \le m \le 0.515$ and $0.9181 \le m \le 0.9187$; if $m \ge 0.9188$, no available switching angles are found.

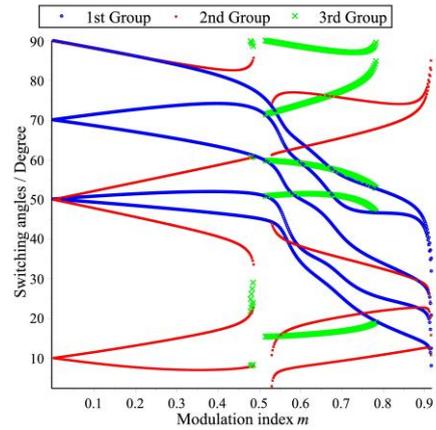

Fig. 3. Switching angle trajectories when $N=5$.

We build a three-level Neutral Point Clamped inverter system using MATLAB/Simulink. The simulation circuit configuration is shown in Fig. 4. The DC side use VSC structure, with parallel connections of capacitance voltage sources in order to maintain the input voltage as a constant. Inputting different switching angles produces corresponding switching signals, and by this way we control the output voltage of the inverter. To verify the effectiveness of the algorithm, no filter is connected to the outlet. A three-phase resistive load is connected to AC side, and the load neutral point is connected to the ground. We measure the voltage waveform at the inverter exit to verify the effectiveness of harmonic suppression. When $m=0.75$, three groups of switching angles: $G_1 = [21.218°, 26.939°, 36.526°, 46.817°, 53.842°]$, $G_2 = [10.055°, 21.255°, 33.889°, 66.911°, 74.966°]$, $G_3 = [17.534°, 49.299°, 54.967°, 79.869°, 87.110°]$ are obtained. We simulate the three groups of switching angles based on the above inverter system. The simulation result of $G_1$ is shown in Fig 5, Fig 6, that of $G_2$ is shown in Fig. 7, Fig. 8, and that of $G_3$ is shown in Fig 9, Fig 10. As can be seen from the Fig. 6, Fig. 8, and Fig. 10, the $5^{th}, 7^{th}, 11^{th}, 13^{th}$ harmonics has been eliminated. This illustrates the validity and effectiveness of this method.

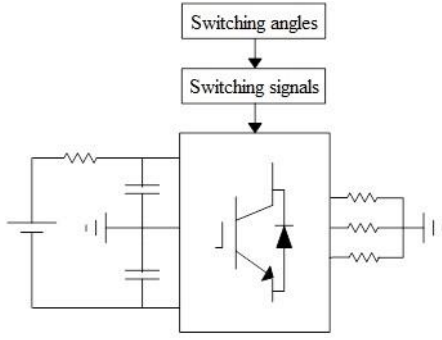

Fig. 4.  Simulation circuit configuration.

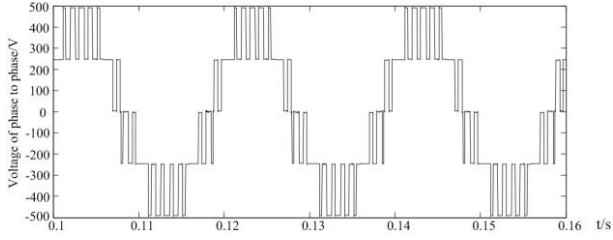

Fig. 5.  Phase to phase output voltage of $G_1$.

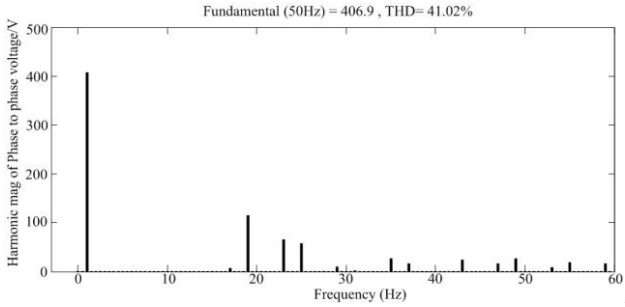

Fig. 6.  Harmonics spectrum of
phase to phase output voltage of $G_1$.

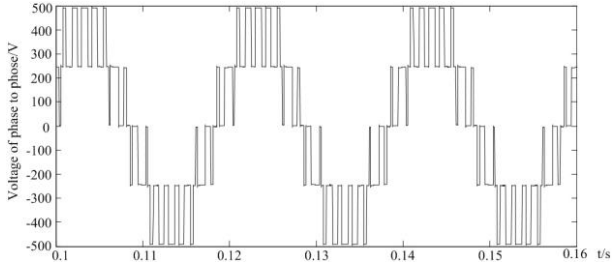

Fig. 7.  Phase to phase output voltage of $G_2$.

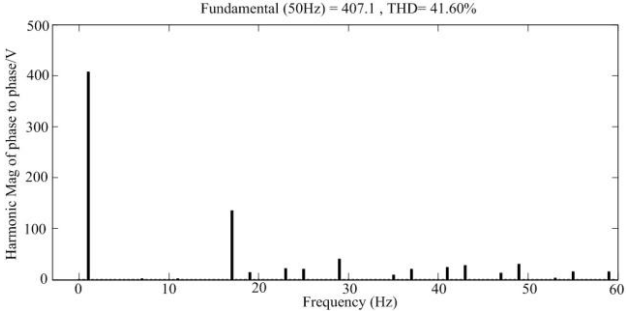

Fig. 8.  Harmonics spectrum of
phase to phase output voltage of $G_2$.

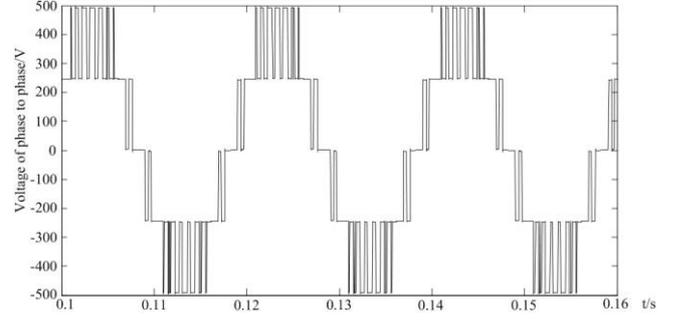

Fig. 9.  Phase to phase output voltage of $G_3$.

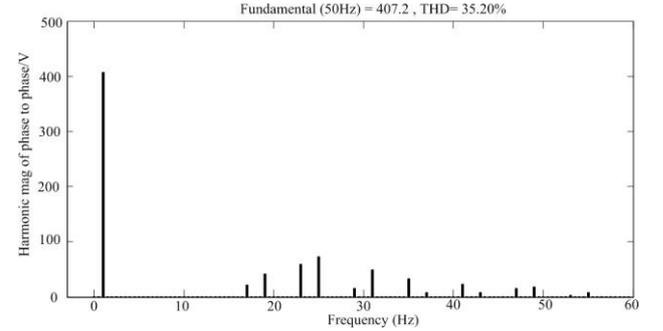

Fig. 10.  Harmonics spectrum of
phase to phase output voltage of $G_3$.

## VI. CONCLUSION

In this paper, based on RUR for solving polynomial systems we presents an algorithm to compute complete solutions to SHEPWM problems. The algorithm can get a "nearly parametric" solution to a SHEPWM problem with the number of switching angles $N$ fixed, and the parameter is the modulation index $m$. Given a specific $m$, only using univariate polynomial solving, we can obtain the complete solutions to the SHEPWM problem. When $m$ changes, the algorithm needn't restart totally, and it speeds up the solving. Experiments show that it's effective to solve SHEPWM problems based on a "nearly parametric" solution.

## APPENDIX

*A. A detailed example with $N = 3$*

The first part of the algorithm in section IV goes as follows.
(1) The transcendental equations are

$$\begin{cases} \dfrac{4}{\pi}\left((\cos\alpha_1 - \cos\alpha_2 + \cos\alpha_3) - m\right) = 0, \\ b_5 = \dfrac{4}{5\pi}\left(\cos 5\alpha_1 - \cos 5\alpha_2 + \cos 5\alpha_3\right) = 0, \\ b_7 = \dfrac{4}{7\pi}\left(\cos 7\alpha_1 - \cos 7\alpha_2 + \cos 7\alpha_3\right) = 0. \end{cases} \quad (8)$$

(2) Let $x_1 = \cos\alpha_1$, $x_2 = -\cos\alpha_2$, $x_3 = \cos\alpha_3$. Translate equations (8) into polynomial equations in $x_1, x_2, x_3$, and get

$$\begin{cases} g_1 = \sum_{i=1}^{3} x_i - m = 0, \\ g_2 = \sum_{i=1}^{3} (5x_i - 20x_i^3 + 16x_i^5) = 0, \\ g_3 = \sum_{i=1}^{3} (-7x_i + 56x_i^3 - 112x_i^5 + 64x_i^7) = 0. \end{cases} \quad (9)$$

(3) Using symmetric polynomial theory, equations (9) are converted into polynomial equations in elementary symmetric polynomials $s_1, s_2, s_3$. Then we substitute $m$ for $s_1$ in $f_2, f_3$, and get

$$\begin{cases} f_1 = s_1 - m = 0, \\ f_2 = 16m^5 - 80m^3 s_2 - 20m^3 + 80m^2 s_3 + 80m s_2^2 \\ \quad + 60m s_2 - 80 s_2 s_3 + 5m - 60 s_3 \\ \quad = 0, \\ f_3 = 64m^7 - 448m^5 s_2 - 112m^5 + 448m^4 s_3 + 896m^3 s_2^2 \\ \quad + 560m^3 s_2 - 1344m^2 s_2 s_3 - 448m s_2^3 + 56m^3 - 560m^2 s_3 \\ \quad - 560m s_2^2 + 448m s_3^2 + 448 s_2^2 s_3 - 168m s_2 + 560 s_2 s_3 \\ \quad - 7m + 168 s_3 \\ \quad = 0. \end{cases} \quad (10)$$

Due to $s_1 - m = 0$, $s_1$ is equivalent to $m$. So actually equations (10) have 1 less variable than equations (9). The number of zeros of equations (9) is 18, and that of equations (10) is 3 (in the algebraic closure of $\mathbb{Q}(m)$). In fact, due to the symmetric property of equations (9), the number of its zeros is 3!=6 times of that of (10) because $N = 3$.

(4) Compute the RUR of $f_2, f_3$ with respect to $m$ and a separating element $t = s_2$. We get $\chi_{t,m}(T)$, $g_{t,m}(1,T)$, $g_{t,m}(s_2,T)$, $g_{t,m}(s_3,T)$, where

$$\begin{cases} \chi_{t,m}(T) = \frac{-3072m^{10} + 15360m^8 - 28160m^6 + 24080m^4 - 9800m^2 + 1575}{35840m^4 - 44800m^2 + 11200} \\ \quad + \frac{1024m^8 - 3904m^6 + 5152m^4 - 2800m^2 + 525}{1792m^4 - 2240m^2 + 560} T \\ \quad + \frac{-576m^6 + 1456m^4 - 1120m^2 + 245}{448m^4 - 560m^2 + 140} T^2 + T^3, \\ g_{t,m}(1,T) = \frac{1024m^8 - 3904m^6 + 5152m^4 - 2800m^2 + 525}{1792m^4 - 2240m^2 + 560} \\ \quad + \frac{-576m^6 + 1456m^4 - 1120m^2 + 245}{224m^4 - 280m^2 + 70} T + 3T^2, \\ g_{t,m}(s_2,T) = \frac{9216m^{10} - 46080m^8 + 84480m^6 - 72240m^4 + 29400m^2 - 4725}{35840m^4 - 44800m^2 + 11200} \\ \quad + \frac{-1024m^8 + 3904m^6 - 5152m^4 + 2800m^2 - 525}{896m^4 - 1120m^2 + 280} T \\ \quad + \frac{576m^6 - 1456m^4 + 1120m^2 - 245}{448m^4 - 560m^2 + 140} T^2, \\ g_{t,m}(s_3,T) = \frac{2048m^{11} - 13056m^9 + 28160m^7 - 27440m^5 - 12600m^3 - 2275}{35840m^4 - 44800m^2 + 11200} \\ \quad + \frac{-2304m^9 + 12480m^7 - 20160m^5 + 12600m^3 - 2625m}{8960m^4 - 11200m^2 + 2800} T \\ \quad + \frac{32m^7 - 140m^4 + 140m^3 - 35}{112m^4 - 140m^2 + 35} T^2, \end{cases} \quad (11)$$

It is a "nearly parametric" solution to the SHEPWM problem.

For a fixed modulation index $m_0$, we can get associated switching angles via univariate polynomial solving. And the second part of the algorithm is illuminated by the following. For example, let $m_0 = 1/2$. (For a specific $m_0$, we always restart the solving process here).

(1) Substitute $m = m_0$ into (11), then we get

$$\begin{cases} \chi_{t,m}^{(m_0)}(T) = 247/2240 + (45/56)T + (47/28)T^2 + T^3, \\ g_{t,m}^{(m_0)}(1,T) = 45/56 + (47/14)T + 3T^2, \\ g_{t,m}^{(m_0)}(s_2,T) = -741/2240 - (45/28)T - (47/28)T^2, \\ g_{t,m}^{(m_0)}(s_3,T) = -449/4480 - (549/1120)T - (33/56)T^2. \end{cases} \quad (12)$$

(2) The three real roots of $\chi_{t,m}^{(m_0)}(T)$ are $T_1 = -0.96307746$, $T_2 = -0.47388307$, $T_3 = -0.24161090$.

(3) Substitute $T = T_i, i = 1,2,3$ into $s_2 = \frac{g_{t,m}^{(m_0)}(s_2,T)}{g_{t,m}^{(m_0)}(1,T)}$,

$s_3 = \frac{g_{t,m}^{(m_0)}(s_3,T)}{g_{t,m}^{(m_0)}(1,T)}$, then the three zeros of equation (10) is got

$s_{1,1} = m_0 = 1/2$, $s_{1,2} = -0.9630774604$, $s_{1,3} = -0.4950353924$;
$s_{2,1} = m_0 = 1/2$, $s_{2,2} = -0.4738830663$, $s_{2,3} = 0.002366801877$;
$s_{3,1} = m_0 = 1/2$, $s_{3,2} = -0.2416109013$, $s_{3,3} = -0.09661712410$.

(13)

According to Theorem 1, we need $s_{i,1} > 0, s_{i,2} < 0, s_{i,3} < 0$, so $s_{2,1}, s_{2,2}, s_{2,3}$ mustn't generate a group of available switching angles.

The roots of $x^3 - s_{1,1}x^2 + s_{1,2}x - s_{1,3} = 0$ are $-0.9859818035$, $0.5194692655$, $0.9665125380$, and there doesn't exist a permutation of the three roots such that $1 > x_{j_1} > -x_{j_2} > x_{j_3} > 0$, so $s_{1,1}, s_{1,2}, s_{1,3}$ mustn't generate a group of available switching angles either.

The roots of $x^3 - s_{3,1}x^2 + s_{3,2}x - s_{3,3} = 0$ are $x_1 = 0.64191435$, $x_2 = -0.465354137$, $x_3 = 0.323439786$, which satisfy $1 > x_1 > -x_2 > x_3 > 0$. So they form a solution of equations (9), then we can get the associated three switching angles $\alpha_1 = 50.06528°$, $\alpha_2 = 62.26686°$, $\alpha_3 = 71.12892°$, which is a solution to equations (8).

From above discussions, we know that when $N = 3$, $m_0 = 1/2$, there only exists one solution to equation (8).

### B. Notions and Conversion relations between symmetric polynomials

**Definition 1** ([18]) A polynomial $f \in k[x_1, \ldots, x_n]$ is *symmetric* if $f(x_{i_1}, \ldots, x_{i_n}) = f(x_1, \ldots, x_n)$ for every possible permutation $x_{i_1}, \ldots, x_{i_n}$ of the variables $x_1, \ldots, x_n$.

**Definition 2** ([18]) Given variables $x_1, \ldots, x_n$, the *elementary symmetric polynomials* $s_1, \ldots, s_n \in k[x_1, \ldots, x_n]$ are defined by the following formulas.

$$\begin{cases} s_1 = x_1 + \cdots + x_n, \\ \vdots \\ s_r = \sum_{i_1 < i_2 < \cdots < i_r} x_{i_1} \cdots x_{i_r}, \\ \vdots \\ s_n = x_1 x_2 \cdots x_n. \end{cases}$$

**Definition 3** ([19]) Given variables $x_1, \ldots, x_n$, the *power sum symmetric polynomials* $p_1, \ldots, p_n \in k[x_1, \ldots, x_n]$ are defined by the following formulas.

$$\begin{cases} s_1 = x_1 + \cdots + x_n, \\ \vdots \\ s_r = x_1^r + \cdots + x_n^r, \\ \vdots \\ s_n = x_1^n + \cdots + x_n^n. \end{cases}$$

The conversion relation from power sum symmetric polynomials to elementary symmetric polynomials can be explicitly represented as follows.[19]

$$\begin{cases} p_1 = s_1, \\ p_2 = s_1^2 - 2s_2, \\ p_3 = s_1^3 - 3s_2 s_1 + 3s_3, \\ p_4 = s_1^4 - 4s_2 s_1^2 + 2s_2^2 - 4s_4, \\ p_5 = s_1^5 - 5s_2 s_1^3 + 5s_3 s_1^2 + 5s_2^2 s_1 - 5s_4 s_1 - 5s_3 s_2 + 5s_5, \\ \vdots \\ p_m = \sum_{r_i=0}^{\lfloor m/i \rfloor} \frac{m(r_1 + r_2 + \cdots + r_n - 1)!}{r_1! r_2! \cdots r_n!} \prod_{i=1}^{n}(-s_i)^{r_i}. \end{cases}$$

### C. Proof of Theorem 1

Proof: Let $s_1^{(k)} = \sum_{i=1}^{k} x_i$, ..., $s_r^{(k)} = \sum_{1 \leq i_1 < \cdots < i_r \leq k} x_{i_1} x_{i_2} \cdots x_{i_r}$, ..., $s_k^{(k)} = x_1 \cdots x_k$, $1 \leq r \leq k$, $1 \leq k \leq n$. If $s_{4d+1}^{(k)} > 0$, $s_{4d+2}^{(k)} < 0$, $s_{4d+3}^{(k)} < 0$, $s_{4d+4}^{(k)} > 0$, $k = 1, \ldots, n$, $d = 0, 1, \ldots$ hold, then it is clear that Theorem 1 holds.

Let $\text{sign}(x) = \begin{cases} 1, & x > 0; \\ 0, & x = 0; \\ -1, & x < 0. \end{cases}$ Apply induction on $r$.

(1) $r = 1$. It is self-evident that $s_1^{(1)}, \ldots, s_1^{(n)} > 0$.

(2) $r = 2$. It can be verified that
$$s_2^{(k)} = x_1 x_2 + \cdots x_1 x_k$$
$$+ x_2 x_3 + \cdots x_2 x_k$$
$$+ \cdots$$
$$+ x_{k-1} x_k$$
$$= x_1(x_2 + \cdots + x_k) + \cdots + x_i(x_{i+1} + \cdots + x_k) + \cdots + x_{k-1} x_k.$$

From $x_1 > -x_2 > \cdots > (-1)^{n+1} x_n > 0$, we know $x_i(x_{i+1} + \cdots + x_k) < 0$, $i = 1, \ldots, k-1$, $2 \leq k \leq n$, hence $s_2^{(k)} < 0$, $k = 2, \ldots, n$.

(3) $r = 3$. It can be verified that
$$s_3^{(k)} = x_1 x_2 (x_3 + \cdots + x_k)$$
$$+ (x_1 + x_2) x_3 (x_4 + \cdots + x_k)$$
$$+ \cdots$$
$$+ (x_1 + \cdots + x_{k-2}) x_{k-1} x_k$$
$$= s_1^{(1)} x_2 (x_3 + \cdots + x_k) + s_1^{(2)} x_3 (x_4 + \cdots + x_k)$$
$$+ \cdots + s_1^{(k-2)} x_{k-1} x_k.$$

Because $s_1^{(k)}$ have the same sign, $\text{sign}(s_3^{(k)}) = \text{sign}(-s_1^{(1)}) < 0$, $k = 3, \ldots, n$.

Generally speaking, we have
$$s_r^{(k)} = s_{r-2}^{(r-2)} x_{r-1}(x_r + \cdots + x_k) + s_{r-2}^{(r-1)} x_r (x_{r+1} + \cdots + x_k)$$
$$+ \cdots + s_{r-2}^{(k-2)} x_{k-1} x_k.$$

Because $s_{r-2}^{(r-2)}, \ldots, s_{r-2}^{(k-2)}$ have the same sign, $\text{sign}(s_r^{(k)}) = \text{sign}(-s_r^{(k-2)}) = \text{sign}(s_r^{(k-4)})$, $k > 4$. In summary, we have $s_{4d+1}^{(k)} > 0, s_{4d+2}^{(k)} < 0, s_{4d+3}^{(k)} < 0, s_{4d+4}^{(k)} > 0, k = 1, \ldots, n, d = 0, 1, \ldots$.